\documentclass[11pt,letterpaper]{amsart}
\usepackage[utf8]{inputenc}
\usepackage{amsmath}
\usepackage{amssymb}
\usepackage{hyperref}
\usepackage{amsthm}
\usepackage{tikz-cd}
\usepackage{cleveref}
\usepackage[shortlabels]{enumitem}
\usepackage{MnSymbol}
\usepackage{pdfpages}
\usepackage{setspace}

\numberwithin{equation}{section}

\newcommand{\img}[2][1]{\begin{gathered}\includegraphics[scale=#1]{#2}\end{gathered}}

\newtheorem{thm}{Theorem}[section]

\newtheorem{ex}[thm]{Example}

\newtheorem{question}[thm]{Question}
\newtheorem{prob}[thm]{Problem}

\newtheorem{thmintro}[thm]{Theorem}
\newtheorem{propintro}[thm]{Proposition}
\newtheorem{conjintro}[thm]{Conjecture}

\newcommand{\bico}{\operatorname{BiCo}}
\newcommand{\Ho}{\operatorname{Ho}}
\newcommand{\mult}{\operatorname{mult}}

\newcommand{\del}{\partial}
\newcommand{\delbar}{\bar{\partial}}

\newcommand{\C}{\mathbb{C}}
\newcommand{\Q}{\mathbb{Q}}
\newcommand{\QQ}{\mathbb{Q}}
\newcommand{\RR}{\mathbb{R}}

\newcommand{\CC}{\mathbb{C}}

\newcommand{\Z}{\mathbb{Z}}
\newcommand{\Id}{\operatorname{Id}}

\newcommand{\Oh}{\mathcal{O}}

\newcommand{\id}{\operatorname{id}}
\newcommand{\cV}{\mathcal{V}}
\newcommand{\cW}{\mathcal{W}}
\newcommand{\cA}{\mathcal{A}}
\newcommand{\im}{\operatorname{im}}

\author{Jonas Stelzig}
\title[Differential forms and invariants of complex manifolds]{Differential forms and invariants of complex manifolds}
\begin{document}

\maketitle
\begin{abstract}
	A survey of some results and open questions related to the following algebraic invariants of compact complex manifolds, that can be obtained from differential forms: cohomology groups, Chern classes, rational homotopy groups, and higher operations.
\end{abstract}
\section*{Introduction}
Compact complex manifolds are ubiquitous in geometry, number theory, and theoretical physics. Yet, there is still a large gap in our understanding of these objects. Broadly speaking, many tools for their study either use only the underlying almost complex structure, or apply only to the special cases of K\"ahler or  projective manifolds. 
In consequence, the geography of general compact complex manifolds, and the boundary regions with the almost complex and K\"ahler worlds, remain largely mysterious.	For example, starting in dimension $3$, not a single example of a compact almost complex manifold which does not admit a complex structure is known -- the most famous instance of this problem being the six-sphere. 
 Similarly, no compact complex manifold of dimension $\geq 3$ is known on which every complex structure is K\"ahler -- this is open even for $\mathbb{CP}^3$.

	A natural first step to amend this state of affairs is to study algebraic invariants that are defined for all compact complex manifolds. On the one hand, one should understand the restrictions imposed upon such invariants by the existence of a K\"ahler metric. On the other hand, one should clarify to what extent they depend on, or can be upgraded to take into account, the presence of an integrable almost complex structure. Our aim here is to give an overview of some results and open questions that are concerned with the following types of such invariants: cohomology groups, Chern classes, rational homotopy groups, and higher operations.
	\[
	\begin{tikzcd}
		\pi_{\geq 2}(X)\otimes\Q&&H_{sing}^{\ast}(X),\, H_{\delbar}^{\ast,\ast}(X),\, H_{BC}^{\ast,\ast}(X),...\\
		&X\ar[ur]\ar[dr]\ar[ul]\ar[dl]&\\
		 \substack{\text{Massey products, }\\A_\infty\text{-structure,...}}\hspace{4ex}&&c_1(X),\, c_2(X),...
	\end{tikzcd}
\]	
	
Our choice of topics is certainly rather biased and incomplete. To name just some of the many omissions, we will never speak about nontrivial fundamental groups and rarely about torsion phenomena, and mostly take coefficients to be $\QQ$, $\RR$, or $\CC$. In this case, all the above invariants can essentially be computed using differential forms, as we will recall.
	
	A guiding question will be that of realization, i.e. which sets of invariants can actually be realized by a space. We will also sometimes briefly touch upon the question of classification, i.e. to what extent do the given invariants determine our space.

	The realization questions we discuss come in two related flavors: a structural one and a numerical one. To fix ideas, let us think about (say, rational) cohomology rings of some class of spaces, a case discussed in more detail below. A structural question is which abstract rings arise as cohomology rings of spaces in the given class. A numerical one is which sequences of numbers arise as the Betti numbers of our class of spaces. The latter is essentially equivalent to asking about all relations that have to hold among the entries of such sequences when they come from a space in our class. In particular, a first approximation would be to ask which linear relations need to hold.

\section{Spaces, manifolds, and almost complex structures} 	
\subsection{Spaces} In the following $X$ will always denote a topological space, which will subsequently be assumed to carry more and more extra structure. Although it is unnecessary in many places, for simplicity of exposition we will generally assume $X$ to be simply connected and of finite type. 

A basic invariant of a topological space is its singular cohomology ring with rational coefficients $H^*_{sing}(X;\QQ)$, or, even more basic, the sequence of Betti numbers $b_k(X):=\dim H^k_{sing}(X;\QQ)$. The answer to the numerical realization question -- which sequences $(b_k)_{k\geq 0}$ of non-negative integers arise as Betti numbers of spaces? -- is clearly: all of them, as one may see by considering wedges of spheres. Structurally, the question becomes: Which graded-commutative graded $\QQ$-algebras arise as the cohomology ring of spaces?

Again the answer is: all of them. To give a precise statement let us introduce some notation. We abbreviate  graded-commutative non-negatively graded $\QQ$-algebra as rational cga. We call a rational cga simply connected if $H^0=\Q$ and $H^1=0$.

\begin{thmintro}\label{realization of coh. rings} For any degree-wise finite-dimensional simply connected rational cga $H^\ast$, there is a simply connected space $X$ with $H_{sing}^\ast(X;\QQ)\cong H^\ast$.
\end{thmintro}
Theorem \ref{realization of coh. rings} can be deduced as a consequence of the deeper results of rational homotopy theory, \cite{Quillen_RHT}, \cite{Sullivan}. Before we state them, let us explain some context.

Recall that singular cohomology is defined to be cohomology of the complex of singular cochains $C^*(X,\QQ)$. The multiplication is defined on the cochain level, but is not commutative there. On the other hand, if $X$ is a manifold, one can associate to it the complex of $\RR$-valued differential forms $A_{X,\RR}^\ast$, which does have a graded-commutative multiplication. It thus has the structure of a graded-commutative differential graded $\RR$-algebra, which will be abbreviated as real cdga. The cohomology of $A_{X,\RR}$, the de Rham cohomology $H_{dR}^\ast(X;\RR):=H(A_{X,\RR},d)$ is thus a real cga and the de Rham theorem gives an isomorphism of cga's: $H^*_{sing}(X;\RR)\cong H^*_{dR}(X;\RR)$.

It is an insight of Sullivan \cite{Sullivan}, based on earlier work of Cartan, Quillen, Thom, Whitney and others \cite{Sullivan_CardeRham}, that one can extend the method of differential forms to rational coefficients and general spaces. Sullivan functorially (cf. \cite{BouGug_PLdR}) attaches a rational cdga of $A_{PL}(X)$ to any space $X$, such that $H^\ast_{sing}(X;\QQ)\cong H^*(A_{PL}(X),d)$. Roughly speaking, $A_{PL}(X)$ consists of polynomial forms on every simplex mapping into $X$. Conversely, there is a `geometric realization' functor, attaching a space $\langle A\rangle$ to any rational cdga $A$. A map of two simply-connected spaces is a rational homotopy equivalence if it induces an isomorphism in cohomology with rational coefficients. Similarly, in the category of rational cdga's one has a natural notion of quasi-isomorphism, namely, a map of cdga's inducing an isomorphism in cohomology. 

Theorem \ref{realization of coh. rings} then follows from the following more general result. 

\begin{thmintro}
	The rational homotopy theory of finite type simply connected spaces is
	equivalent to the homotopy theory of finite type cohomologically simply
	connected rational cdga's.
	
%
\end{thmintro}
In a more technical language: Sullivan's $A_{PL}$-forms and the geometric realization induce an equivalence of categories between the full subcategories of the category of spaces localized at rational homotopy equivalences, resp. the category of rational cdga's localized at quasi-isomorphisms, which are spanned by those objects that are of finite type and (cohomologically) simply connected.

In view of this theorem, it becomes a natural question to find an optimal representative for the quasi-isomorphism class of a rational cdga $A$. We recall two ways of doing so.

One way is by Sullivan's minimal models: for any connected cdga $A$, there exists a minimal cofibrant model, i.e. a quasi-isomorphism from another cdga $M_A\to A$ such that $M_A=\Lambda V$ and $V$ has a well-ordered basis $(v_i)$ with $v_i<v_j$ for  $\deg v_i <\deg v_j$ such that $d v_i$ is a linear combination of products of two or more lower-order generators. Every self homotopy equivalence of such an $M_A$ is an actual isomorphism and so the homotopy class of $A$ can be identified with the isomorphism class of $M_A$.

Another way is via higher operations and homotopy transfer: To motivate this, recall the classical triple Massey products \cite{Massey_mexico}, \cite{Massey_holn}. For any three pure degree classes $[a],~[b],~[c]\in H^\ast(A)$ s.t. $ab=dx$, $bc=dy$ one obtains a new class $\langle a,b,c\rangle:=[ay-(-1)^{|a|}xc]\in H^{\ast}(A)$, well-defined up to the ideal generated by $[a]$ and $[c]$ and depending only on the cohomology classes of $a$, $b$, $c$. This can be considered as a partially defined ternary operation $H^\ast(A)^{\otimes 3}\dashedrightarrow H^\ast(A)$. The theory of homotopy transfer gives a slightly more involved way of equipping $H(A)$ with everywhere defined $n$-ary operations $m_k$ for $k\geq 1$ (with $m_1=0$ and $m_2$ being the usual product), satisfying certain compatibility conditions making it into a so-called $C_\infty$-algebra. Any cdga with $m_1=d$, $m_2$ the product and all higher operations trivial is in particular a $C_\infty$-algebra. Then, the homotopy class of $A$ can be identified with the isomorphism class of the $C_\infty$-algebra $(H^\ast(A),0,m_2,m_3,m_4,...)$.
Let us summarize this discussion in the following theorem:
\begin{thmintro}\label{thm: mimos}
Sullivan's $A_{PL}$-functor, building a minimal model, resp. homotopy transfer, induce bijections preserving the cohomology ring between the following sets of (quasi-)isomorphism classes, where we restrict to (cohomologically) simply connected objects of finite type in each case:
\tikzcdset{scale cd/.style={every label/.append style={scale=#1},
		cells={nodes={scale=#1}}}}
\[
\begin{tikzcd}[scale cd=0.7]
	&\left\{\begin{tabular}{c}
		Rational homotopy equivalence \\
		classes of spaces
	\end{tabular}\right\}\ar[d,leftrightarrow]&\\
&\left\{\begin{tabular}{c}
	Quasi-isomorphism classes\\
	of rational cdga's\\
\end{tabular}\right\}\ar[ld,leftrightarrow]\ar[rd,leftrightarrow]&\\
\left\{\begin{tabular}{c}
	Isomorphism classes \\
	of minimal cdga's\\
\end{tabular}\right\}
&&\left\{\begin{tabular}{c}
	Isomorphism classes \\
	of $C_\infty$-algebras \\
	with $m_1=0$ 
\end{tabular}\right\}\,.
\end{tikzcd}
\]
\end{thmintro}

\subsection{Compact manifolds}
Let us now consider compact smooth manifold of dimension $n$, without boundary. Again one can ask the realization questions: Which collections of non-negative integers are Betti numbers of compact smooth $n$-folds? Which rational cga's are cohomology algebras of compact smooth $n$-folds?

One obvious obstruction comes from Poincar\'e duality: For any compact smooth $n$-fold, the $n$-dimensional cohomology is one-dimensional $H^n(X;\Q)\cong \QQ$ and multiplication yields a perfect pairing $H^k(X;\QQ)\times H^{n-k}(X,\QQ)\to H^n(X;\Q)\cong \Q$. Any finite dimensional cga satisfying these properties will be called a (rational) PD-algebra of formal dimension $n$. For such algebras, the linear relation $b_k=b_{n-k}$ needs to hold. Furthermore, in even dimensions $n=2k$, the induced pairing in middle cohomology is $(-1)^k$ symmetric. In particular, on all PD-algebras of formal dimension $n$ with $n\equiv 2~\operatorname{mod}~(4)$, the congruence $b_{n/2}\equiv 0~\operatorname{mod}~(2)$ holds.

Not every collection of non-negative numbers $b_0,...,b_n$ with $b_0=1$ satisfying the above duality constraints can be obtained as the Betti numbers of a connected compact $n$-fold as we will see shortly. However, we have a positive answer to the linearized numerical realization problem:

\begin{thmintro}[\cite{KoSchr_HodgeKaehl}]
	The relations $b_k=b_{n-k}$ and $b_{n/2}\equiv 0~\operatorname{mod}~(2)$ if $n\equiv 2~\operatorname{mod}~(4)$ are the only relations that hold universally for all Betti numbers of compact smooth $n$-manifolds.
\end{thmintro}


If one wants to realize not just Betti numbers, but cohomology rings, there are, in dimensions $n=4k$, more subtle restrictions coming from Poincar\'e duality and characteristic classes. First, since Poincar\'e duality holds integrally, the non-degenerate symmetric pairing in middle degree needs to be rationally equivalent to one of the form $\sum \pm y_i^2$.

Next, recall that for any (stably) complex vector bundle $\cV$ on $X$, one has a classifying map $f_\cV:X\to BU$ and the cohomology of $BU$ is a polynomial algebra on generators $c_i$ situated in degrees $2i$. The Chern classes of the vector bundle $\cV$ are then the classes $c_i(\cV):=f_{\cV}^*c_i\in H^{2i}(X;\Z)$. For a real vector bundle $\cW$, one defines its Pontryagin classes by $p_i(\cW):=(-1)^ic_{2i}(\cW\otimes\C)\in H^{4i}(X;\Z)$. In particular, for any smooth manifold $X$, one writes $p_i(X):=p_i(TX)$. If $X$ is compact with fundamental class $[X]$ and of dimension $n$ divisible by $4$, one thus obtains for every partition $\tau$ of $\frac n4$, say $\frac n 4 =\sum_{i=1}^{\frac n 4} k_i\cdot i$, the Pontryagin number 
\[p_{\tau}(X):=\left\langle p_1(X)^{k_1}\cdot...\cdot p_{n/4}(X)^{k_{n/4}},[X] \right\rangle \in\Z\,.\]

By Hirzebruch's signature theorem, in every dimension divisible by $4$, the signature $\sigma$ of the middle degree pairing can be computed as a universal linear expression with rational coefficients in the Pontryagin numbers. For example, dropping evaluation on the fundamental class in the notation (taken from \cite{OEIS_LPoly}):
\[{\renewcommand{\arraystretch}{1.4}
\begin{tabular}{l|c}
$n$&$\sigma$\\\hline
$4$&$\frac 13 p_1$\\
$8$&$\frac{1}{45}(7 p_2-p_1^2)$\\
$12$&$\frac{1}{945} (62 p_3-13 p_2 p_1+2 p_1^3)$\\
$16$&$\frac{1}{14175} (381 p_4-71 p_3 p_1-19 p_2^2+22 p_2 p_1^2-3 p_1^4)$\\
$20$&$\frac{1}{467775} (5110 p_5-919 p_4 p_1-336 p_3 p_2+237 p_3 p_1^2+127 p_2^2 p_1-83 p_2 p_1^3+10 p_1^5)$.
\end{tabular}}
\]
This shows that there are necessary congruences for the Pontryagin classes and hence not all sequences of Betti numbers satisfying the necessary linear relations are realized by manifolds. For example:\footnote{This argument is taken from Hirzebruch's 1953 manuscript `The index of an oriented manifold and the Todd genus of an almost complex manifold', see \cite{Hirz_GesAbh}.} If there were a compact $12$-fold with Betti numbers $b_0=b_6=b_{12}=1$ and $0$ else, it would follow that $\frac{62}{945}\in\Z$, which is absurd.
With this in mind, one has the following version of Theorem \ref{realization of coh. rings} for compact smooth manifolds:

\begin{thmintro}[Sullivan-Barge realization theorem \cite{Sullivan}, \cite{Barge}]\,\\
Let $H^\ast$ be a PD-algebra of formal dimension $n\geq 5$ and $H^1=0$. Fix some classes $p_i\in H^{4i}$.
\begin{enumerate}
	\item If $n\not\equiv 0 ~\operatorname{mod}~(4)$, in every rational homotopy type with cohomology ring $H^\ast$, there exists a  compact smooth manifold $X$ with $H(X)\cong H^\ast$ and $p_i(X)=p_i$.
	\item If $n\equiv 0 ~\operatorname{mod}~(4)$, the same conclusion holds if we assume furthermore that there is a choice of fundamental class in $(H^n)^\vee$ such that the pairing induced on $H^{n/2}$ is rationally equivalent to one of the form $\sum \pm y_i^2$ and that the signature can be computed from Hirzebruch's expressions evaluated in the `Pontryagin numbers' formed formally from $p_i$.
\end{enumerate}

\end{thmintro}

In particular, in all dimensions $\geq 5$, whether a rational homotopy type contains a simply connected smooth closed manifold only depends on the rational cohomology ring.

In fact, one can even solve the classification problem up to finite ambiguity \cite{Sullivan}: Roughly speaking, the statement is that if one takes into account integral information, one may refine the above to a statement roughly saying that the map
\[
\{\text{closed manifolds}\mid  \pi_1=0\}\longrightarrow\{\text{minimal model, lattices, torsion, } p_i\}/_\sim
\]
induces, in dimension $\geq 5$, a finite-to-one map on diffeomorphism classes.

\subsection{Almost complex manifolds}
A manifold $X$ together with an endomorphism $J$ of the tangent bundle squaring to $-\Id$ is called an almost complex manifold. By definition, the tangent bundle is a complex vector bundle, and so the real dimension of the manifold is even, say $2n$. Furthermore, one has canonically defined Chern classes $c_i(X):=c_i(X,J):=c_i(TX)\in H^{2i}_{sing}(X;\Z)$. If $X$ is in addition compact, one can thus form Chern numbers
\[
c_\tau(X):=\left\langle c_1(X)^{k_1}\cdot ... \cdot c_{n}(X)^{k_{n}}, [X]\right\rangle\in\Z
\]
for any partition $\tau$ of $n$, written as $ n  = k_1+k_2\cdot 2 + ... + k_{n}\cdot n$.

Again these numbers have to satisfy certain universal congruences on all almost complex manifolds of a given dimension,  characterizing the image of the map $\Omega^U\to H_\ast (BU,\Q)$, resp. in the case $c_1=0$ the image of $\Omega^{SU}\to H_\ast (BSU,\Q)$. Roughly speaking, they are all integrality conditions that come out of the Atiyah-Singer index theorem. A precise statement is as follows, we call these the Stong congruences following \cite{Mili}:

\begin{thmintro}[Stong congruences, \cite{AtiyahHirzebruch_CohOpsChaChla}, \cite{Stong_CN_I}, \cite{Stong_CN_II}, \cite{Hattori_CN}]\,\\
	Let $X$ be a closed almost complex manifold of dimension $2n$.
	\begin{enumerate}
		\item For any $n$, the numbers
		\[\langle z\cdot Td(X),[X]\rangle\]
		are integers for every polynomial $z\in\Z[e_1,e_2,...]$, where $e_i$ are the elementary symmetric polynomials in the variables $e^{x_j}-1$, where $x_j$ are given by formally writing $1+c_1+c_2+...=\Pi_j(1+x_j)$ and $Td(X)$ denotes the Todd polynomial evaluated on $c_1,c_2,...$.
		\item If $n\equiv 0 ~ (4)$ and $c_1=0$, furthermore the numbers
		\[
		\frac 1 2\langle z\cdot \hat{A}(X),[X]\rangle
		\]
		are integers for every $z\in\Z[e_1^p,e_2^p,...]$, where $e_i^p$ are the elementary symmetric polynomials in the variables $e^{x_j}+e^{-x_j}-2$, where the $x_j$ are given by formally writing $1+p_1+p_2+....=\Pi_j(1+x_j^2)$ and $\hat{A}(X)$ denotes the $\hat{A}$ polynomial evaluated on $p_1,p_2,...$.
	\end{enumerate}
\end{thmintro}

The first set of congruences in low dimensions are as follows, \cite{Hirz_KoMa}, \cite{Mili}:
\[{\renewcommand{\arraystretch}{1.4}
	\begin{tabular}{l|c}
		$2n$&universal congruences\\\hline
		$2$&$c_1\equiv 0~\operatorname{mod}~(2)$\\
		$4$&$c_2+c_1^2\equiv 0~\operatorname{mod}~(12)$\\
		$6$&$c_3\equiv c_1^3\equiv0~\operatorname{mod}~ (2), \quad c_2c_1\equiv 0~\operatorname{mod}~(24)$\\
		$8$&$-c_4+c_3c_1+3c_2^2+4c_2c_1^2-c_1^4\equiv 0~\operatorname{mod}~(720)$,\\
		&$c_2c_1^2+2c_1^4\equiv 0 ~\operatorname{mod}~(12),\quad -2c_4+c_3c_1\equiv 0~\operatorname{mod}~(4)$\\
		$10$&$c_5+c_4c_1\equiv 0 ~\operatorname{mod}~(12),\quad 9c_5+c_4c_1+8c_3c_1^2+4c_2c_1^3\equiv 0~\operatorname{mod}~ (24)$,\\
		&$-8c_4c_1+8c_3c_1^2+12c_2c_1^2-5c_2c_1^3+15c_1^5\equiv 0~\operatorname{mod}~(24)$,\\
		&$6c_3c_1^2+c_2c_1^3+c_1^5\equiv 0~\operatorname{mod}~(12),$\\
		&$ -c_4c_1+c_3c_1^2+3c_2^2c_1-c_2c_1^3\equiv 0~ \operatorname{mod}~(1440)$
\end{tabular}}
\]
Then one has the following almost complex version of the Sullivan-Barge theorem, which, after the formulation and proof had been left to the reader in \cite{Sullivan}, has been worked out fully in \cite{Mili}:

\begin{thmintro}[Almost complex realization, \cite{Mili}]\label{thm acSB}
	Let $A$ be a rational cdga of finite type with $H^1=0$ and satisfying rational Poincar\'e duality with formal dimension $2n$. Fix a choice of a class $0\neq[X]\in H_{2n}(A)=H^{2n}(A)^\vee$ and classes $c_i\in H^{2i}(X,\Q)$. Then the homotopy type of $A$ contains a closed, simply connected, almost complex manifold $X$ realizing $[X]$ as fundamental class and the $c_i$ as Chern classes if and only if the following conditions are satisfied:
	\begin{enumerate}
		\item The `Chern numbers' formed formally with the classes $c_i$ and $[X]$ satisfy the Stong congruences.
		\item If $n$ is odd, the `intersection form' on $H^{n}(A,\Q)$ is rationally equivalent to one of the form $\sum \pm y_i^2$ and the signature can be computed from the Hirzebruch $L$-polynomial evaluated in the `Pontryagin numbers' formed formally from $c_i$ and $[X]$. 
		\item The Euler characteristic equals the top Chern number: $\langle c_{n},[X]\rangle =\sum (-1)^k b_k(A)$.
	\end{enumerate}
\end{thmintro}
We note that, as in the smooth case, the realizability of a simply connected rational homotopy type depends only on the cohomology ring. Unlike the smooth case, the above data combined with integral information do not finitely determine the almost complex manifold up to pseudoholomorphic equivalence. For instance, the action of the diffeomorphism group of a given real manifold of real dimension $\geq 4$ on the space of almost complex structures has an infinite-dimensional orbit space.

\section{Complex manifolds}

We will now restrict our focus to integrable almost complex structures on compact manifolds. While it may still be unreasonable to expect a finite-to-one classification in general, the situation improves in at least two ways: First, the moduli space of integrable almost complex structures is locally finite dimensional\footnote{This is not true globally, however, see e.g. \cite{Brieskorn_PnBundleP1}.}  \cite{KodairaSpencer_DefComI}, \cite{KodariaSpencer_DefComII}, \cite{KodairaNirenbergSpencer}. Second, there are many more finite-dimensional cohomology theories, which depend on the complex structure and capture some information on the moduli.\footnote{The question of generalizing the definitions of the cohomologies we discuss below to general (compact) almost complex manifolds is currently actively being pursued, to some extent prompted by \cite[Problem 20]{Hirzebruch}. Without meaning to give a full survey, there are roughly two directions: (1) Definitions as spaces of harmonic forms for certain Laplace operators. The dimension of the so-defined vector spaces are finite (on compact manifolds), but generally depend on the choice of a metric, with  some notable exceptions in dimension $4$, see e.g. \cite{HZ}, \cite{HZ_AK}, \cite{AK}, \cite{PioTo_BCAH}, \cite{Holt_BCddb4}. (2) Metric independent definitions, see e.g. \cite{LZ_tc}, \cite{CiWi_HdR4}, \cite{Stelzigandfriends}, \cite{CaGuGu_traDol},  \cite{SilTo_BCA-acs}. These tend to not always be finite dimensional, even on compact manifolds \cite{Stelzigandfriends}.} We survey these next. 

\subsection{A zoo of cohomology theories}\label{sec: zoo} On any almost complex manifold $(X,J)$, the cdga of $\C$-differential forms $A_{X}=A_{X,\RR}\otimes\C$ carries a bigrading $A_X:=\bigoplus A_X^{p,q}$. It is induced by the splitting of the bundle of $\C$-valued $1$-forms $\cA_X^1=\cA_X^{1,0}\oplus \cA_X^{0,1}$ into $i$- and $-i$-eigenbundles for the endomorphism $J$. With respect to this bigrading, the differential splits into components $d=\del + \delbar$ of bidegree $(1,0)$ and $(0,1)$, if and only if $J$ is integrable. We restrict to this case from now on. By the equation $d^2=0$, one has $\del^2=\delbar^2=\del\delbar+\delbar\del=0$. In other words, the complex of forms on a complex manifold is the total complex underlying a bicomplex. The multiplication is compatible with the bigrading, so that $A_X$ carries the structure of a graded-commutative bidifferential, bigraded algebra, which we will abbreviate as cbba. The antilinear conjugation action on $A_X$ interchanges $A_X^{p,q}$ with $A_X^{q,p}$ and $\del$ with $\delbar$. A cbba with such an antilinear isomorphism will be called an $\RR$-cbba. A smooth map $f:X\to Y$ of complex manifolds is holomorphic iff it respects the bigrading, i.e. $f^*(A_Y^{p,q})\subseteq A_X^{p,q}$.

From $A_X$ one can build various holomorphic invariants of $X$. We write down the following definitions  only for $A_X$, but they are meaningful for any bicomplex (or $\RR$-cbba), and we will use this in the next sections.

The most well-known holomorphic cohomological invariant is perhaps Dolbeault cohomology \cite{Dolbeault_cohomologie}, controlling existence and uniqueness of the $\delbar$-equation $x=\delbar y$. It is defined as the graded-commutative bigraded algebra, obtained as the cohomology of the columns of the bicomplex $(A_X,\del,\delbar)$:
\[
H_{\delbar}(X):=\frac{\ker\delbar}{\im\delbar}\,.
\]
It can be identified with the sheaf cohomology of the sheaf of holomorphic functions $H_{\delbar}^{p,q}(X)=H^q(X,\Omega^p)$. 

Next, note that on the de Rham cohomology with complex coefficients $H_{dR}(X):=H_{dR}(X;\C)$, one has a multiplicative filtration, induced from the column filtration on $A_X$, i.e.
\[
F^pH_{dR}^\ast(X):=\im \left(\bigoplus_{r\geq p} A^{r,s}_X\cap\ker d\longrightarrow H_{dR}^\ast(X)\right),
\]
and another filtration $\bar F$, computed analogously from the rows of $A_X$, which is conjugate to $F$. These filtrations generally depend on the complex structure.

As for any bicomplex, there is a spectral sequence, the Fr\"olicher spectral sequence \cite{Fro},
\[
E_1^{p,q}(X)=H_{\delbar}^{p,q}(X)\Longrightarrow (H_{dR}^{p+q}(X),F)\,,
\]
which does not degenerate in general and so the higher pages $E_r^{p,q}(X)$ give additional invariants of the bihomolomorphism type of $X$.

On any compact $X$ the vector spaces $H^{p,q}_{\delbar}(X)$ (and hence all later pages) are finite dimensional and satisfy Serre duality, i.e. for any connected compact $X$ of dimension $n$, one has $H^{n,n}_{\delbar}(X)\cong \C$ and the multiplication
\[
H^{p,q}_{\delbar}(X)\times H^{n-p,n-q}_{\delbar}(X)\longrightarrow H^{n,n}_{\delbar}(X)\cong\C
\]
is a perfect pairing. The vector spaces $H_{\delbar}(X)$ are generally not conjugation invariant, rather, conjugation swaps $H_{\delbar}^{p,q}(X)$ with $H_{\del}^{q,p}(X)$, where the latter vector space is defined analogously but exchanging $\delbar$ by $\del$.

Another natural bigraded cohomology algebra, which is conjugation invariant, is the Bott-Chern cohomology
\[
H_{BC}(X):=\frac{\ker\del\cap\ker\delbar}{\im\del\delbar}\,.
\]
It was defined by Bott and Chern in \cite{BottChern} studying generalizations of Nevanlinna theory and has found many other uses since.\footnote{To highlight a few, the secondary characteristic classes of \cite{BottChern}, have proven useful in the study of Hermite Einstein metrics \cite{Don_BC} but also Arakelov theory and refined versions of the Riemann-Roch theorem \cite{GiSou_CCI}, \cite{GiSou_CCII}, \cite{Gi_ICM}, \cite{Bismut_RR}. Bott-Chern cohomology is sometimes a more natural setup to generalize results from the K\"ahler case to arbitrary compact complex manifolds. For example, the statement and proof of the Castelnuovo de Franchis theorem as e.g. in \cite[Ch. IV, §5]{BaHuPeVdV_CCS} generalize if one replaces Dolbeault cohomology by Bott-Chern cohomology. Much more deeply, there is a hermitian analogue of the Calabi-conjecture, see e.g. \cite{TosattiWeinkove_cplxMA}, \cite{SzkelyhidiTosattiWeinkove_Gaud}; one obtains an interesting variant of Calabi-Yau manifolds via the condition $c_{1,BC}(TX)=0$ \cite{Tosatti_NKCY}.} It is an initial cohomology in the sense that it has natural maps to $H_{\delbar}, H_{\del}, H_{dR}$ and all other cohomology theories discussed later on. For any holomorphic vector bundle $\cV$, the Chern classes can naturally be lifted to Bott-Chern cohomology in the following sense: There are classes $c_{i,BC}(\cV)\in H_{BC}^{i,i}(X)$, compatible with pullback by holomorphic maps, such that under the natural map $H_{BC}(X)\to H_{dR}(X)$, they coincide with the image of $c_i(\cV)$ under the map $H_{sing}(X;\Z)\to H_{dR}(X)$.

Bott-Chern cohomology generally does not satisfy a Serre-type duality. Rather, it pairs non-degenerately with Aeppli cohomology \cite{Aeppli_exact}, \cite{Aeppli_coh}
\[
H_A(X):=\frac{\ker\del\delbar}{\im\del+\im\delbar}\,.
\]
These vector spaces are also stable under conjugation and are a `terminal' cohomology in the sense that they receive natural maps from $H_{\delbar}, H_{\del}, H_{dR}$ and all other cohomology theories discussed later on. 

Bott-Chern and Aeppli cohomology groups arise as hypercohomology groups in certain degrees of the complexes of sheaves \cite{Schw_BC}, \cite{Demailly_book}, \cite{Bigolin_Aepp}, \cite{Bigolin_Oss}
\[
\mathcal{S}_{p,q}=( \Oh_X+\bar\Oh_X\to \Omega^1_X\oplus\bar\Omega^1_X\to...\to \Omega_X^{p-1}\oplus \bar\Omega_X^{p-1}\to\bar\Omega_X^{p}\to ...\to\bar\Omega_X^{q-1})
\]  
which are closely related to the complexes computing the Deligne cohomology groups \cite{Beil_DelCoh}. Other cohomology groups of these complexes appear for example in the classification of holomorphic string algebroids \cite{GarRubTip_holstringalgebroids}, or holomorphic higher operations \cite{MilSt_bigrform}. Some general properties of these cohomologies are established in \cite{SteSC} and \cite{Piov_ABC}.

One can continue this list with more cohomology theories, arising naturally from certain geometric contexts. E.g. if one extends $J$ as an algebra-automorphism to the vector space $A_X$ of all forms and sets $d^c:=J^{-1}dJ$, one may obtain new complexes $(\ker d^c, d)$ and $(A/\im d^c, d)$ and corresponding singly graded cohomology theories 
\[
H_{\ker d^c}^\ast (X):=H^\ast(\ker d^c, d), \quad H_{A/\im d^c}^\ast(X):= H^\ast(A/\im d^c, d),
\]
which are also finite-dimensional and in duality, \cite{SteWi}. Then there are the groups defined by Varouchas \cite{Var_Prop}, and higher-page analogues of those and Aeppli and Bott-Chern cohomology \cite{PoSteU}, \cite{PoSteUb}, \cite{PoSteUc} and possibly more which we have not mentioned.

All of the above cohomologies are finite dimensional, and so in addition to the Betti and Hodge numbers $b_k(X):=\dim H_{dR}^k(X)$ and $h_{\delbar}^{p,q}(X):=\dim H_{\delbar}^{p,q}(X)$, one obtains a vast collection of numerical invariants from the dimensions of the respective cohomologies: $h_{BC}^{p,q}:=\dim H_{BC}^{p,q}(X)$, $h_A^{p,q}:=H_A^{p,q}(X)$ etc. One also has a bigraded refinement of the Betti numbers using the filtrations, setting $b_k^{p,q}(X):=\dim \operatorname{gr}_F^p\operatorname{gr}_{\bar F}^qH_{dR}^k(X)$.

These numbers satisfy certain universal relations. E.g. on all compact complex manifolds of a given dimension $n$, one has linear relations induced by duality and the real structure, such as
\[
e_{r,\delbar}^{p,q}=e_{r,\delbar}^{n-p,n-q}=e_{r,\del}^{q,p}, \qquad h_{BC}^{p,q}=h_{BC}^{q,p}=h_A^{n-p,n-q},\qquad b_k^{p,q}=b_k^{q,p}=b_{2n-k}^{n-p,n-q}\,,
\]
but also other linear relations such as
\[
h^{0,n}_{\delbar} = h_{BC}^{0,n},\qquad h^{0,n-1}_{\delbar} = h^{0,n-1}_A,\qquad h^{n-1,0}_{\delbar} = h^{n-1,0}_{BC}\,.
\]
There is then the following natural question, which we will come back to later:

\begin{question}\label{que: lin rel coh inv}
	What are the linear relations that hold universally between the dimensions of the various cohomologies associated to compact complex manifolds of a given dimension?
\end{question}

We note that in dimension $n=2$, one also has a non-linear polynomial relation
\[
(h^{0,1}_{\delbar} - h_{BC}^{0,1})(h^{0,1}_{\delbar}- h_{BC}^{0,1}-1)=0,
\]
but it is unknown whether universal non-linear polynomial relations exist in higher dimensions.\footnote{Relatedly, I do not know if there is a universal upper bound of $h^{1,0}$ by an expression involving $h^{0,1}$ in any given dimension. E.g. on curves and surfaces $h^{1,0}\leq h^{0,1}$, with equality always holding in the K\"ahler case. The Iwasawa manifold gives an example in complex dimension $3$ for which $h^{1,0}=h^{0,1}+1$, but I am not aware of an example in dimension $3$ where the difference $h^{1,0}-h^{0,1}$ can be bigger. It can be arbitrarily small, as one can show the calculation of the Hodge numbers of a twistor space $Z:=Z(M)$ of a self-dual compact Riemannian $4$-fold $M$ in \cite{EaSi}. In fact, these yield $h^{1,0}(Z)=0$, $h^{0,1}(Z)=b_1(M)$).}

Further, there are universally valid inequalities. Let us denote the total dimension of any of the above cohomology theories by a letter without sub- or superscripts indicating a degree, so $b:=\sum_{k=0}^{2n} b_k$, $h_{BC}:=\sum_{p,q=0}^n h_{BC}^{p,q}$, $e_i:=\sum_{p,q=0}^{n}\dim E_i^{p,q}$, etc. 

\begin{thmintro}[\cite{AnTo}, \cite{PoSteUb}, \cite{SteWi}] The following two sets of inequalities hold on all compact complex manifolds:
\begin{equation}\label{eqn: bc-dc-delb-b}
h_{BC}\geq h_{\ker d^c}\geq h_{\delbar}\geq b,
\end{equation}
and, for any fixed $r\geq 0$,
\begin{equation}\label{eqn: BC-e-b}
h_{BC}\geq \sum_{i=1}^{r} e_i(X)-(r-1)b\,.
\end{equation}
\end{thmintro}
It is an interesting question to find manifolds on which equality holds. E.g. one has $h_{\ker d^c}=h_{\delbar}$ for Vaisman manifolds, all complex surfaces and complex parallelizable manifolds with solvable Lie algebra of holomorphic vector fields \cite{PoSteUb}, \cite{KSt_parallel}, \cite{SteWi}. In the first two cases, in addition $h_{\delbar}=b$, while in the last case, in addition $h_{BC}=h_{\ker d^c}$. There are also simply connected examples for both cases \cite{KSt_nK}. There is, however, by no means any known kind of classification of manifolds satisfying a particular case of equality. In particular, it is unknown whether there are any manifolds satisfying equality in \ref{eqn: BC-e-b} for $r\geq 2$, but not for $r-1$.

\subsection{The compact K\"ahler case}

An important class of compact manifolds are those admitting a K\"ahler metric. This includes all complex submanifolds of projective space. 

In this case, the cohomological story greatly simplifies: Roughly speaking, all cohomologies are determined by Dolbeault cohomology. More precisely:

\begin{propintro}[The $\del\delbar$-Lemma, \cite{DGMS}]
	For any bicomplex $A=(A,\del,\delbar)$, the following assertions are equivalent:
	\begin{enumerate}
		\item For any $a\in A$ such that $\del a = \delbar a=0$ and $a=db$ for some $b\in A$, there exists a $c\in A$ such that $a=\del\delbar c$.
		\item There is an isomorphism $A\cong A^{sq}\oplus A^{dot}$, where $A^{sq}$ is a direct sum of squares, i.e. bicomplexes of the form
		\begin{equation}\label{eqn: square}
		\begin{tikzcd}
			\C\ar[r]&\C\\
			\C\ar[u]\ar[r]&\C,\ar[u]
		\end{tikzcd}
		\end{equation}
		where all arrows are $\pm\id$ and all other maps vanish, and $A^{dot}$ is a direct sum of dots, i.e. one-dimensional bicomplexes with all differentials being zero.
		\item All maps in the diagram 
		\begin{equation}\label{eqn: diagram}
			\begin{tikzcd}
				&H_{BC}(A)\ar[ld]\ar[d] \ar[rd]&\\
				H_{\delbar}(A)\ar[rd]\ar[r,Rightarrow]&H_{dR}(A)\ar[d]&H_{\del}(A)\ar[ld]\ar[l,Rightarrow]\\
				&H_A(A)&
		\end{tikzcd}\end{equation}
	are isomorphisms.
	\item The spectral sequences in the previous diagram degenerate, and the filtrations on de Rham cohomology are $n$-opposed, i.e. $b_k^{p,q}=0$ unless $k=p+q$.
	\end{enumerate}
	Moreover, for a compact K\"ahler manifold $X$, the bicomplex $A_X$ satisfies these conditions.
\end{propintro} 
A complex manifold $X$ for which $A_X$ satisfies the above conditions is called a $\del\delbar$-manifold. Compact K\"ahler manifolds are $\del\delbar$-manifolds, but the converse it not true. A broader class is for example given by those manifolds bimeromorphic to compact K\"ahler manifolds, i.e. Fujiki's class $\mathcal{C}$, but also these do not exhaust all $\del\delbar$-manifolds, see e.g. \cite{Friedm_ddbarClemens}, \cite{Li_PolHodgeClemens}, \cite{KSt_nK}.

One readily checks that all cohomologies introduced above are compatible with direct sums, vanish on squares and are one-dimensional on dots. Thus, as soon as a cohomology allows to reconstruct the information on the position of the dots, e.g. Dolbeault cohomology, or de Rham cohomology with its filtrations, it determines all other ones on $\del\delbar$-manifolds. 
In such cases, the answer to Question \ref{que: lin rel coh inv} is given by the following theorem:

\begin{thmintro}[\cite{KoSchr_HodgeKaehl}]\label{thm: KS HC rels}
	The only universal $\QQ$-linear relations between Hodge, Betti and Chern numbers on compact K\"ahler manifolds of dimension $n$ are:
	\begin{enumerate}
		\item (Real structure and duality) $h^{p,q}=h^{q,p}=h^{n-p, n-q}$
		\item (Hodge decomposition) $b_k=\sum_{p+q=k}h^{p,q}$
		\item (Hirzebruch Riemann Roch) $\chi_p=Td_p$
	\end{enumerate}
	Here $\chi_p:=\sum_p (-1)^qh^{p,q}$ and $Td_p$ is the $p$-th Todd-genus, a certain linear combination of Chern numbers. 
\end{thmintro}

It is stated here for compact K\"ahler manifolds as in \cite{KoSchr_HodgeKaehl}, but from the proof in \cite{KoSchr_HodgeKaehl} it is clear that it remains true when considering the smaller (resp. larger) classes of projective, resp. $\del\delbar$-manifolds. The result (for Hodge numbers) was extended to polynomial relations in Paulsen and Schreieder \cite{PaSchr_Hodge}. See also \cite{vDdB_HodgeRing-p}, \cite{vDdBPa_Hodge-p} for characteristic $p$ versions of these results.

\subsection{Hirzebruch's question}\label{sec: Hirz}
Given a biholomorphism $X\to Y$ of complex manifolds, it induces an isomorphism $H(Y)\to H(X)$ for $H$ any of the cohomologies introduced above. In particular, the cohomology dimensions are invariants of the biholomorphism class of a manifold.

On the other hand, some linear combinations of these cohomology dimensions are even topological invariants. This is most obvious for $h_{\delbar}^{0,0}$ which counts connected components. By a spectral sequence argument, one sees that one can compute the Euler characteristic as $\chi=\sum_{p,q}(-1)^{p+q}h^{p,q}_{\delbar}$. More deeply, by the Hodge index theorem, for general compact complex manifolds $\sigma=\sum_{p,q}(-1)^q h^{p,q}_{\delbar}(X)$ holds 
 and so the combination of Hodge numbers on the right is invariant under orientation preserving homeomorphisms. In view of such relations, Hirzebruch asked the following question in 1954 \cite{Hirzebruch}:
\[
\img[0.55]{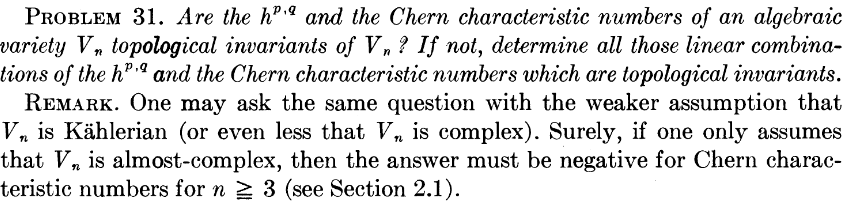}
\]
There are several remarks in order: First, one may understand `topological invariant' in at least four a priori different ways, namely invariant under homeomorphisms or diffeomorphisms of the underlying manifold which may or may not be required to preserve the orientation.
Next, as mentioned before there are universal linear relations between the Hodge and Chern numbers of compact complex manifolds in a given dimension, described for compact K\"ahler manifolds in Theorem \ref{thm: KS HC rels} and one should therefore answer the question modulo these relations.

Some extreme cases of the problem were quickly resolved: E.g. in 1958, Hirzebruch and Borel exhibited an example of diffeomorphic $5$-folds with distinct $c_1^5$ \cite{BorHirz_charchlahomsp_I-III}, giving a negative answer to the initial question. A determination of topologically invariant Chern numbers of mere almost complex manifolds was given in \cite{Kahn_Chernhtt}.

Much more recently, the problem concerning the Chern numbers only was solved by Kotschick \cite{Ko_ChaNu}, \cite{KotProb}, following earlier work \cite{Kot_Orientrevers}, \cite{Kot_OrientationsGeometrisations} and finally the original problem for both Hodge and Chern numbers of projective varieties was solved by Kotschick and Schreieder:

\begin{thmintro}[\cite{KoSchr_HodgeKaehl}]
	A rational linear combination of Hodge and Chern numbers of compact projective manifolds of dimension $n\geq 3$ is
	\begin{enumerate}
		\item an oriented homeomorphism or diffeomorphism invariant if and only if it reduces to a linear combination of the Betti and Pontryagin numbers modulo the universal relations of Theorem \ref{thm: KS HC rels}.
		\item an unoriented homeomorphism or diffeomorphism invariant if and only if it reduces to a linear combination of the Betti numbers modulo the universal relations of Theorem \ref{thm: KS HC rels}.
	\end{enumerate}
\end{thmintro}

It is clear from the proof of this theorem that one can replace projective manifolds by K\"ahler manifolds or even by $\del\delbar$-manifolds and the statement remains unaffected. 

It is proved in \cite{SteLR} that in the above formulation the answer remains the same when relaxing to all compact complex manifolds of a given dimension $n\geq 3$, only that the universal relations are different (e.g. $h^{p,q}_{\delbar}\neq h^{q,p}_{\delbar}$ in general).\footnote{This makes the proof nontrivial, as there are many linear combinations which are equivalent on K\"ahler manifolds, but not on general compact complex manifolds, thus one has to construct more examples showing that certain combinations are not topological invariants.} We will furthermore see that there are more general questions to be asked in the general compact complex realm, taking into account all cohomological invariants instead of just the Hodge numbers.


\subsection{Common motives} \label{sec: motives}
All the previously discussed cohomologies behave similarly in many geometric situations, e.g. they can be computed as one would expect from the de Rham case for projective bundles and blow-ups \cite{St_blowup}. However, they do not all satisfy a straightforward K\"unneth formula, see e.g. \cite{ChiRas_AstKaehlABC}, \cite{SteLR}, \cite{StePHT}.

Naturally, one expects a more fundamental invariant lurking behind all these cohomology groups. In some sense there is an obvious answer: They are all computed from the bicomplex of forms $(A_X,\del,\delbar)$. In fact, one readily checks that they are well-defined for a general bicomplex $(A,\del,\delbar)$ and that they are naturally compatible with direct sums $H(\bigoplus A_i)\cong \bigoplus H(A_i)$. Finally, they all vanish on squares as in Diagram \eqref{eqn: square}.

These observations may be seen as a motivation for the following definitions: We denote by $\bico$ the category of all bicomplexes and by $\Ho(\bico)$ the quotient category, where all morphisms that factor over a direct sum of squares are set to zero. This is called the homotopy category, or derived category, of bicomplexes. An additive functor $H:\bico\to \operatorname{Ad}$ to some additive category $\operatorname{Ad}$ is called a cohomological functor if it factors through this derived category. In other words, it has to commute with finite direct sums and vanish on direct sums of squares. A functor from complex manifolds to $\operatorname{Ad}$ is called cohomological if it factors, via $X\to A_X$ through a cohomological functor. Examples of such functors are of course all the above cohomology theories and also diagrams of those, such as Diagram \eqref{eqn: diagram}.

Now, the values of cohomological functors on objects $A$ are determined by the isomorphism type of $A$ in $\Ho(\bico)$, and so instead of studying the cohomology theories of a complex manifold one at a time one should determine this isomorphism type for the bicomplex $(A_X,\del,\delbar)$, i.e. determine $A_X$ `up to squares'. To do this, one needs a good understanding of what maps of bicomplexes induce isomorphisms in $\Ho(\bico)$ and how to represent an isomorphism class in $\Ho(\bico)$, so let us discuss these points next.

We call a map of bicomplexes inducing an isomorphism in $\Ho(\bico)$ a bigraded quasi-isomorphism. By definition it has a quasi-inverse, i.e. a map in the other direction such that the compositions both ways differ from the identity by maps factoring through direct sums of squares. This is the bicomplex analogue of a map of complexes invertible up to chain homotopy. One may characterize these maps cohomologically as follows:

\begin{propintro}[\cite{StStrDbl}, \cite{StePHT}]Let $f:A\to B$ a map of bicomplexes. 
	\begin{enumerate}
		\item The map $f$ is a bigraded quasi-isomorphism if and only if the induced maps $H_{BC}(A)\to H_{BC}(B)$ and $H_A(A)\to H_A(B)$ are isomorphisms.
		\item Assume that for every fixed $k$, $H_{BC}^{p,q}(A)$ and $H_{BC}^{p,q}(B)$ 
		are nonzero for only finitely many bidegrees with $p+q=k$. Then, the map $f$ is a bigraded quasi-isomorphism if and only if the induced maps $H_{\delbar}(A)\to H_{\delbar}(B)$ and $H_{\del}(A)\to H_{\del}(B)$ are isomorphisms.
	\end{enumerate}
\end{propintro}

Thus, a map of bounded bicomplexes is a bigraded quasi-isomorphism if and only if it is a quasi-isomorphism in row and column cohomology and, using the real structure, a map of complex manifolds $f:X\to Y$ induces a bigraded quasi-isomorphism $A_Y\to A_X$ if and only if it induces an isomorphism in Dolbeault cohomology. For both $A$ and $B$ satisfying the $\del\delbar$-Lemma, a bigraded quasi-isomorphism is the same as a map of bicomplexes which is also a usual quasi-isomorphism (i.e. it induces an isomorphism in total cohomology). In view of the characterization in terms of Bott-Chern and Aeppli cohomology, which measure the existence and uniqueness of solutions to the $\partial\bar\partial$-equation $x=\partial\bar\partial y$, and to avoid confusion with other possible meanings of `bigraded quasi-isomorphism', we will sometimes also use the name pluripotential quasi-isomorphism to denote the same concept.

To describe the pluripotential quasi-isomorphism type, or even the actual isomorphism type, of a bicomplex, the following theorem is useful:

\begin{thmintro}[\cite{KhQi}, \cite{StStrDbl}]\label{thm: Dec}
	Every bicomplex is a direct sum of indecomposable bicomplexes. Every indecomposable bicomplex is isomorphic to a square or a zigzag. 
\end{thmintro}
Here, a zigzag is a bicomplex concentrated in at most two antidiagonals, s.t. all nonzero vector spaces are $\C$, all nonzero maps are the identity and all nonzero vector spaces are connected by a chain of maps. E.g.:
\begin{equation}\label{eqn: zigzags}
\begin{tikzcd}
\C\ar[r]&\C&\\
&\C\ar[r]\ar[u]&\C
\end{tikzcd}\qquad\text{ or }\qquad 
\begin{tikzcd}
\C\ar[r]&\C\\
&\C\ar[u]
\end{tikzcd}
\end{equation}
The total dimension of a zigzag will be called its length (so the above have length $4$ and $3$).\footnote{This maybe counterintuitive convention, counting components and not arrows for the length, has the advantage that the tensor product of two odd zigzags is (up to squares) again odd, while the product of any zigzag with an even zigzag is (up to squares) a direct sum of even zigzgas, cf. \cite{StStrDbl}.} Zigzags of length one are called dots and those of length two are called lines. For every length $\geq 2$, there are, up to shift, two zigzags of that length, and then there are zigzags which may be infinite in one, or both directions.

Given some bicomplex $A$ and a zigzag $Z$, denote by $\mult_Z(A)$ the number of direct summands isomorphic to $Z$ in a decomposition as in \Cref{thm: Dec}. This number does not depend on the chosen decomposition. Thus, any bicomplex $A$ can be written as a direct sum of squares and zigzags $A\cong A^{sq}\oplus A^{zig}$. The bigraded quasi-isomorphism type (knowing $A$ up to squares) is uniquely determined by the isomorphism type of $A^{zig}$ which is again uniquely determined by the collection of numbers $\mult_Z(A)$ for all zigzags $A$. One can encode this information in a `checkerboard' diagram. E.g. suppose $A$ is the direct sum of the two zigzags in Figure \eqref{eqn: zigzags}, where we assume the top left corner of the first to sit in degree $(0,1)$ and the top left corner of the second to sit in degree $(3,1)$. Then we may depict this as
\[
A\simeq \img[1.2]{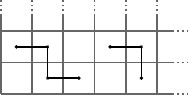}
\]
Now the situation discussed in the previous section has become much more transparent: The value (on objects) of a cohomological functor $H$ is determined by its value on all zigzags. E.g. we urge the reader to stop for a moment and convince themselves that the Fr\"olicher spectral sequence(s) degenerates on the first page for every odd zigzag, while for any even zigzag of length $2r$, the de Rham cohomology vanishes and there is a nonzero differential on page $r$ of the row- or column Fr\"olicher spectral sequence. In particular, as a consequence of finite-dimensionality of Dolbeault cohomology, the multiplicities of all zigzags (and hence the dimensions of all other cohomological functors discussed) are finite on compact complex manifolds.

Properties of the quasi-isomorphism type of $A_X$ will have shadows in the various cohomological functors. For instance real structure and duality hold on the bicomplex level: For any bicomplex $(A,\del_1,\del_2)$ denote by $\bar{A}$ the conjugate bicomplex, which has the same underlying vector space and total differential, but twisted multiplication and bigrading, i.e. $\bar{A}^{p,q}=A^{q,p}$ with conjugate $\C$-linear structure.
Then, there is an isomorphism
\begin{equation}\label{eqn: real structure}
A_X\cong \bar{A}_X.	
\end{equation}

Similarly, denote by $DA[n]$ the $n$-th dual bicomplex, i.e. $(DA[n])^{p,q}=Hom(A^{n-p,n-q},\C)$ with total differential given (up to sign) by precomposition with the differential of $A$. This is a module over $A$ by precomposition with the multiplication. For any compact complex $n$-fold $X$, integration over $X$ yields a canonical closed element $\int_X\in DA_X[n]^{0,0}$ and the composition with the module structure yields a bigraded quasi-isomorphism
\begin{equation}\label{eqn: duality}
A_X\simeq DA_X[n].
\end{equation}
Thus, for any cohomological functor one has canonical isomorphisms \[H(A_X)\cong H(\bar{A}_X)\cong H(DA_X[n]).
\] For each particular $H$, it is an easy exercise to work out an expression for the two terms on the right and one may thus recover all the known results about real structures and duality for the individual cohomological functors.

Any numerical invariant of compact complex manifolds of the form $h=\dim\circ H$ for some cohomological functor $H$, is a sum of the invariants $\mult_Z(\_)$ by additivity of cohomological functors. Consequently, all universal relations between all such numerical invariants $h$ on compact complex manifolds of a given dimension, it is necessary and sufficient to determine all universal relations between the numbers $\mult_Z(\_)$ in a given dimension. The following are all known relations, which follow from classical results (in particular, Serre duality and the cohomological properties of compact complex surfaces). See \cite{StStrDbl} for a treatment in this language and references.

\begin{thmintro}\label{thm: known rels}
	Let $X$ be a compact complex $n$-fold. Then the following universal linear relations hold:
\begin{enumerate}
	\item[(R1)]\label{R1} (Real structure) Let $\sigma$ be the involution which associates to every zigzag its mirror along the diagonal $p=q$. Then $\mult_Z(A_X)=\mult_{\sigma Z}(A_X)$ for all zigzags $Z$.
	\item[(R2)]\label{R2} (Duality) Let $\tau$ be the involution which associates to every zigzag its mirror along the antidiagonal $p+q=n$. Then $\mult_Z(A_X)=\mult_{\tau Z}(A_X)$ for all zigzags $Z$.
	\item[(R3)]\label{R3} (Only dots in the corners) Let $Z$ denote any zigzag of length $\geq 2$ which has a nonzero component in degree $(0,0)$, $(n,n)$, $(n,0)$ or $(0,n)$. Then $\mult_Z(A_X)=0$.
	\item[(R4)]\label{R4} (Fr\"olicher degeneration in dimension $2$) For $n=2$, and $Z$ any zigzag of length $2$, $\mult_Z(A_X)=0$.
	\item[(R5)]\label{R5} For $n=2$, if $Z$ denotes the following zigzag, with top left vector space sitting in degree $(0,1)$:
	\[
	\begin{tikzcd}
		\C\ar[r]&\C\\
		&\C,\ar[u]
	\end{tikzcd}
	\] one has the quadratic relation
	 \[
	\mult_Z(A_X)\cdot (\mult_Z(A_X)-1)=0,
	\]	
	or, equivalently, $\mult_Z(A_X)\in \{0,1\}$.
\end{enumerate}
\end{thmintro}
The following more precise version of Question \ref{que: lin rel coh inv} is open, but partial results are obtained in \cite{SteLR}:
\begin{question}\label{que: rels zigzs}
	Are all universal rational linear (resp. polynomial) relations between cohomological invariants of compact complex manifolds of a given dimension a consequence of the relations \textit{(R1)--(R4)} (resp. \textit{(R1)--(R5)})?
\end{question}

Using the first two relations above, one may refine the `checkerboard' notation for compact complex manifolds of a given dimension. Namely, for a complex $n$-dimensional manifold consider an $(n+1)\times(n+1)$ board and instead of single zigzag write the entire $\Z/2\times\Z/2=\langle \tau,\sigma\rangle$ orbit on one board. E.g. in this notation, the bicomplex of forms of a connected compact complex curve $\Sigma_g$ of genus $g$ looks as follows
\[
A_{\Sigma_g}\simeq\img{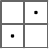}~\oplus~ \img{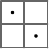}^{~\oplus g},
\]
and that of a connected compact complex surface $X$ looks as follows:
\[
A_X\simeq \img{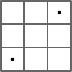}~\oplus~\img{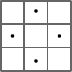}^{~\oplus \frac{b_1-\varepsilon}{2}}~\oplus~ \img{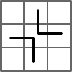}^{~\oplus\varepsilon}~\oplus~\img{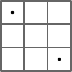}^{~\oplus b_2^++\varepsilon}~\oplus~\img{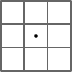}^{~\oplus b_2^-+1-\varepsilon}~,
\]
where $\sigma=b_2^+-b_2^-$ denotes the signature and 
$\varepsilon$ is zero for $b_1(X)$ even (K\"ahler case) and $1$ otherwise.

To give a positive answer to Question \ref{que: rels zigzs}, it essentially remains to solve the following construction problems (cf. \cite{SteLR}):

\begin{prob}\label{prob: hahs}
	For even $n\geq 4$, construct an $n$-dimensional compact complex manifold $X_n$ with $b_{n-1}^{n-1,n-1}(X_n)=1$, i.e. supporting a nonzero $n-1$ de Rham class, unique up to scalar, which can be represented by both a holomorphic and an antiholomorphic form.
\end{prob}

I.e. one is looking for a manifold with an indecomposable summand of the following form in the bicomplex of differential forms:
\[\img{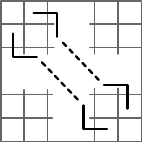}
\]
For odd $n\geq 3$, these exist.
\begin{prob}\label{prob: extremal FSS diffs}
	For every $n\geq 3$, construct a $n$-dimensional compact complex manifold $X_n$ with nonvanishing differentials on page $E_{n-1}$ starting in degree $(0,n-1)$ or $(0,n-2)$.
\end{prob}

I.e. one needs to construct manifolds realizing the following direct summands in their bicomplexes, where for readability we omit zigzags determined via the real structure: 
\[
\img{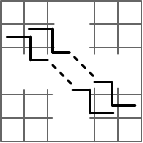} \quad \text{ and }\quad \img{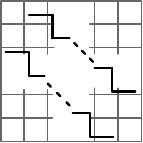}~.
\]
For $n=3$, the first case is known to exist and the second one corresponds to a (hypothetical) threefold with a nontrivial differential $E^{0,1}_2\to E_2^{2,0}$. 

\subsection{Topological invariants}

One notes that the examples in the last section imply that for curves and surfaces, the entire bigraded quasi-isomorphism type is determined by the (oriented) topological manifold underlying $X$. This is far from true in dimensions $\geq 3$. For instance if one takes $X=N/\Gamma$ the quotient of an even-dimensional nilpotent Lie-group modulo a lattice, there are generally many left-invariant complex structures which have distinct cohomological invariants like Hodge numbers, and hence distinct bigraded quasi-isomorphism type, but are all deformation equivalent. LeBrun's examples \cite{LeBrun_TopChern3}, \cite{Mili_HodgevsChern} show there can even be infinitely complex structures with pairwise distinct Hodge numbers on the same smooth manifold.


This, together with Hirzebruch's question begs the following more general question:

\begin{question}
	Which linear combinations of multiplicities of zigzags (and Chern numbers) are topological invariants of compact complex manifolds in dimension $\geq 3$?
\end{question}

E.g. the Betti numbers are a sum of the multiplicities of certain odd zigzags. On the other hand, not all multiplicities of odd zigzags are topological invariants as one may see e.g. in small deformations of the Iwasawa manifold \cite{Nak_par}, \cite{Ang_Iwa}. In analogy with what happens in the Hodge-case, one may conjecture:

\begin{conjintro}\label{conj}
Modulo universal relations, a linear combination of zigzag multiplicities and Chern numbers of compact complex manifolds in dimension $n\geq 3$ is an (orientation preserving) homeomorphism invariant iff it is an (orientation preserving) diffeomorphism invariant iff it is a linear combination of the Betti numbers (and the Pontryagin numbers).
\end{conjintro}

\subsection{Bimeromorphism invariants}

Apart from topological equivalence (i.e. homeomorphism or diffeomorphism, with or without preserving the orientation), there is another natural notion of equivalence between complex manifold, namely that of bimeromorphic equivalence. The typical example of a bimeromorphic map is that of a blow-up $\widetilde{X}\to X$ in some submanifold $Z\subset X$. 
By the deep results of \cite{Wlod}, \cite{AKMW}, two manifolds are related by a bimeromorphic map iff they are related by a chain of roofs of blow-ups in smooth centers. Thus, to show that some number, or isomorphism class of algebraic object, is invariant under bimeromorphisms, it is necessary and sufficient to show it is does not change under blow-ups. The bigraded quasi-isomorphism class (and hence every cohomology) of the blow-up $\widetilde{X}$ of a smooth complex submanifold $Z\subseteq X$ of codimension $r\geq 2$ can be computed as follows \cite{St_blowup}:
\[
A_{\widetilde{X}}\simeq A_X\oplus \bigoplus_{i=1}^{r-1} A_Z[i].
\]
Together with relations \textit{(R1)--(R5)}, this implies:
\begin{thmintro}[\cite{St_blowup}, \cite{StStrDbl}]
	Let $n\geq 2$. The multiplicities of the following zigzags are bimeromorphic invariants:
	\begin{enumerate}
		\item All zigzags having a nonzero component in the boundary degrees $\{0,n\}\times \{0,...,n\}\cup \{0,...,n\}\times \{0,n\}$
		\[
		\img[0.5]{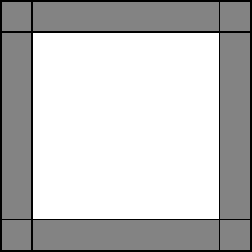}
		\]
		\item All zigzags of length $\geq 2$ which have a component in the `secondary corners' $\{(1,1), (n-1,n-1), (1,n-1), (n-1,1)\}$
		\[
		\img[0.5]{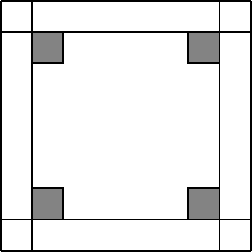}
		\]
		\item If $n\leq 4$, the multiplicities of all even length zigzags, i.e. in addition to the above, the multiplicities of
		\[
		\img[0.8]{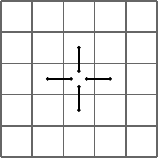}
		\]
		are bimeromorphic invariants.
		\end{enumerate}
\end{thmintro}
The pictures given here might be helpful when reading \cite[§10]{SteLR}. In particular, this theorem recovers the classical results of the bimeromorphic invariance of the Hodge numbers $h^{p,0}$ and $h^{0,q}$, but also exhibits further invariants away from the boundary. Naturally, one may then ask:

\begin{question}\label{que: bimero zigs}
	Is any linear combination of cohomological invariants which is a bimeromorphic invariant of compact complex manifolds of a given dimension a linear combination of the above multiplicities of zigzags?
\end{question}

Again, this question has to be read modulo the (still not completely determined) universal relations. 
The analogous question for Hodge numbers alone (and also including Chern numbers) has a positive answer as was confirmed by Kotschick and Schreieder \cite{KoSchr_HodgeKaehl} in the K\"ahler case and for general compact complex manifolds in \cite{SteLR}. What is missing for a general resolution is an answer to Question \ref{que: rels zigzs}.

\subsection{Multiplicative matters}

In the previous discussion we focused on numerical invariants. There are now many structural realization questions one could ask. For instance, which bigraded rings arise as Dolbeault cohomology or Bott-Chern cohomology rings, which diagrams of the form \eqref{eqn: diagram} arise from compact complex manifolds? Which rational homotopy types contain compact complex manifolds? 

If one asks for realization by compact K\"ahler or projective manifolds, on the one hand the whole cohomological situation is determined by Dolbeault cohomology, on the other hand there is further structure, e.g. coming from the Hard Lefschetz theorem. Known necessary conditions for positive answer to the realizability of cohomology rings have been discussed for example in \cite{Vois_HodgeAlg}, \cite{Vois_Kod}, \cite{Vois_hotKod}.

In the case of general compact complex manifolds, one again faces the situation that there are many distinct cohomologies and corresponding variants of realization questions one could consider. For instance, analogues of rational homotopy theory that replace ordinary cohomology by a particular complex cohomology have to some extent been developed in the Dolbeault or Bott-Chern setting, see e.g. \cite{NT}, \cite{Tan_ModDolFib}, \cite{HaTa} and \cite{AnToBCform}. 

As in the additive case, the need for a universal invariant arises.
A by now obvious candidate is the bicomplex of forms together with its multiplication, turning it into a graded-commutative, bigraded, bidifferential algebra with real structure ($\RR$-cbba) and satisfying Serre duality, up to (multiplicative, real) bigraded quasi-isomorphism.

\begin{question}\label{que: complex realization}
	Consider a collection $(A_\Q,A_\C,h,\{c_{i}, c_{i,BC}\})$ where
	\begin{enumerate}
		\item $A_\Q$ is a rational PD-cdga of formal dimension $2n\geq 6$, and $c_{i}\in H^{2i}(A)$, such that the conditions for almost complex realization in Theorem \ref{thm acSB} are satisfied.
		\item $A_\C$ is an $\RR$-cbba satisfying Serre duality of formal dimension $n$, i.e. \eqref{eqn: real structure} and \eqref{eqn: duality} hold. The underlying bicomplex of $A$ further satisfies \textit{(R3)} of Theorem \ref{thm: known rels}.\footnote{In fact, one should require slightly more: For any compact complex $n$-manifold $X$, the intersection pairing makes $H^{n,0}_{BC}(X)\oplus H^{0,n}_{BC}(X)$ into a polarized (real) Hodge structure. This implies the $(n,0)$ part of \textit{(R3)}.}
		\item $h$ is a chain of quasi-isomorphisms of cdga's between $A_\Q\otimes\C$ and $A_\C$, which is compatible with the real structures 
		on both sides.
		\item $c_{i,BC}\in H_{BC}^{i,i}(A_\C)$ are classes whose images in $H_{dR}(A_\C)$ are identified with $c_{i}$ via the isomorphism $H_{dR}(A_\Q)\otimes\C\cong H_{dR}(A_\C)$.
	\end{enumerate}
	Is there a compact complex $n$-manifold $X$ such that Sullivan's PL-forms, the smooth $\C$-valued forms, the zigzag of quasi-isomorphism given by the PL de Rham theorem, and the Chern classes in rational and Bott-Chern cohomology realize these data, up to an appropriate notion of quasi-isomorphism? 
	 If so, to what extent is $X$ determined by these data?
\end{question}

Much is unknown about the above question. For instance, just as it is not known whether any almost complex manifold in dimension $\geq 6$ is complex, it is not known whether one could always `complete' only the data $(A_\Q, \{c_{i,dR}\})$ to a tuple realizable by a complex manifold. In other words:

\begin{question}\label{que: cga-cbba}
	Consider a rational PD-cdga $A$ which satisfies all the condition for almost complex realization. Is it quasi-isomorphic (over $\RR$) to an $\RR$-cbba satisfying the conditions of Question \ref{que: complex realization} $(2)$, in particular Serre duality?
\end{question}

One may view the question as stated above as an algebraic, or rational homotopy theoretic, version of the question of topological obstructions to the existence of complex structures. Without duality assumption, a positive answer is given in \cite{StePHT}. On the other hand, the stricter question that asks for the existence of an $\RR$-cbba such that the underlying bicomplex satisfies some additive conditions, generally will have a negative answer. The prototypical result in this direction is the famous
\begin{thmintro}[\cite{DGMS}, \cite{Sullivan}]
	Any $\partial\bar\partial$-manifold is rationally formal, i.e. there is a chain of quasi-isomorphisms of cdga's between $A_{PL}(X)$ and $H^\ast_{sing}(X;\QQ)$, where the latter is considered as a cdga with trivial differential.
\end{thmintro}
 A famous proof of this theorem from \cite{DGMS} proceeds as follows: By the $\partial\bar\partial$-Lemma, both maps in diagram of cdga's
 \begin{equation}\label{diag: formality DGMS}
 \begin{tikzcd}
 	&(\ker d^c,d)\ar{ld}\ar{rd}&\\
 	(A_X,d)&&(H_{d^c}(X),0)
 \end{tikzcd}
 \end{equation}
 are quasi-isomorphisms (e.g. one may check this on every indecomposable summand). Since this connects the forms to an algebra with trivial differential, the proof is complete.\footnote{Alternatively, note that $J$ yields a canonical isomorphism $H_{d^c}(X)\cong H_d(X)$.} Note that this argument uses only the $\partial\bar\partial$-property and not that the cbba comes from geometry. Thus a cdga can only be quasi-isomorphic to an $\RR$-cbba satisfying the $\del\delbar$-Lemma, if it is formal.
 
	Inspired by this result, one may ask if there are more general additive conditions than the $\del\delbar$-Lemma on the bicomplex
which are incompatible with certain multiplicative types of cbba's. This is indeed the case as we show in \cite{SteWi}. Instead of reproducing the general statement, let us give two examples from \cite{SteWi} to illustrate the types of questions one can answer with it: 
 
 \begin{ex}
 Consider a filiform nilmanifold $M=G/\Gamma$ where $\Gamma$ is a lattice in the simply connected Lie group $G$ associated with the cdga of left-invariant forms
 \[\Lambda(\eta^1,...,\eta^6)\quad d\eta^1=d\eta^2=0,~d\eta^k=\eta^1\eta^{k-1}\text{ for }k=3,...,6.
 \]
 Like any even-dimensional nilmanifold, $M$ admits an almost complex structure (e.g. put $J\eta^{2k}=\eta^{2k-1}$). It is known that $M$ does not admit left-invariant complex structures \cite{GM}, and it is unknown whether it admits any complex structures.  But what if we impose conditions on the pluripotential quasi-isomorphism type for the bicomplex of forms $(A_M, \del ,\delbar)$? For instance, is there a complex structure with the  following bicomplex possible?
 \[
 A_M\simeq\img{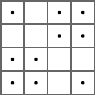}\oplus \img{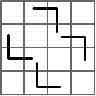}\oplus  \img{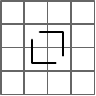}
 \]
 Note that this would yield the correct Betti numbers and have a pure Hodge structure on $H^1$. 
\end{ex}

\begin{ex}
	Let $N=G/\Gamma$ be a nilmanifold with structure equations 
	\begin{align*}
		d\eta^3 &= \eta^1 \eta^2  &  d\eta^4 &= \eta^1 \eta^3 \\
		d\eta^5 &= \eta^2 \eta^3  &  d\eta^6 &= \eta^1 \eta^4 + \eta^2  \eta^5.
	\end{align*}
	Any such nilmanifold has a left invariant complex structure, cf. \cite{Sal}. According to \cite{COUV} (p. 4, Theorem $2.1$) there are two left-invariant complex structures on $N$. In fact, one may compute that for each of them, the bicomplex looks as follows:
	\[
	A_N\simeq\img{heresy}\oplus \img{S311_3}\oplus\img{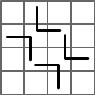}\oplus \img{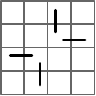}
	\]
	The last summand means that there is a differential on page $1$ of the Fr\"olicher spectral sequence for these structures. One may ask whether there is any complex structure (not necessary left invariant) on $N$ with the same summands except the last. This would have the same bi-filtered de Rham cohomology, but degenerate Fr\"olicher spectral sequence.
\end{ex}
 It turns out, the answer to the questions in both examples is no in a  rather strong sense: 
 \begin{thm}[\cite{SteWi}]
 Structures as asked for in the previous two examples cannot exist on any manifold in the real homotopy type of $M$, resp. $N$.
 \end{thm} 
 The techniques in proving this are not limited to the real homotopy types of nilmanifolds, but also apply to suitable highly connected homotopy types for example.

 Let us survey the main idea behind this type of result: The starting point is the observation is that the diagram \eqref{diag: formality DGMS} exists on any complex manifold. In general, the maps will not be quasi-isomorphisms, but the induced differential on $H_{d^c}(X)$ will still be trivial and one can make a few statements about the two induced maps $H(\ker d^c)\to H(X)$, e.g. they will have the same rank, as one may verify by checking on each indecomposable summand. The next observation is that if one replaces $A_X$ by a quasi-isomorphic algebra $A$, there exists again a diagram
\[ \begin{tikzcd}
 	&B\ar{ld}\ar{rd}&\\
A&&H_{d^c}(A)
 \end{tikzcd}
\]
where the induced maps in cohomology have exactly the same ranks as in \eqref{diag: formality DGMS}. Finally, what exactly the ranks are is determined by which indecomposable summands occur in $A_X$. Now, not every cdga can support such a diagram with arbitrary ranks of the induced maps. Intuitively, if there is such a diagram with the maps close to being isomorphisms, there cannot be many Massey products.\\

\subsection{A pluripotential version of rational homotopy theory} In analogy with Theorem \ref{thm: mimos}, in \cite{StePHT}, canonical representatives for the pluripotential quasi-isomorphism class of a cohomologically simply connected $\RR$-cbba, meaning $H^0_{BC}=H^0_{A}=\C$ and $H^1_{BC}=H^1_A=0$, are found:

\begin{thmintro}[\cite{StePHT}]\label{thm: big mimos}
	The natural forgetful map from left to right yields a canonical bijection
	\[
	\left\{	\begin{tabular}{c}
	Isomorphism classes of\\
	simply connected\\
	minimal $\RR$-cbba's
\end{tabular}\right\}\cong	\left\{	\begin{tabular}{c}
Bigraded quasi-isomorphism classes of\\
cohomologically simply connected\\
$\RR$-cbba's
\end{tabular}\right\}\,.
	\]
\end{thmintro}
The pluripotential variant of homotopy transfer and the $C_\infty$-side of Theorem \ref{thm: mimos} is currently being developed by Anna Sopena-Gilboy (forthcoming, see also \cite{CiGaSo}).

Moreover, in \cite{StePHT}, a model category structure on the category of cbba's is exhibited, for which the forgetful functor to cdga's respects weak equivalences and (co-)fibrations. These results in particular answers positively a question of Sullivan of whether it is possible to build models in the rational homotopy theory sense for complex manifolds which are compatible with bigrading and real structure. In fact, it does so with respect to the strong notion of pluripotential quasi-isomorphism (as opposed to de Rham quasi-isomorphism).

 Further, using these structured models, for any compact simply connected complex manifold, one may define groups that relate to the complexified duals of the homotopy groups just as the various complex cohomology groups relate to de Rham cohomology, e.g.

\begin{equation}\label{diag: complex homotopy}
		\begin{tikzcd}
		&\pi_{BC}^{p,q}(X,x)\ar[ld]\ar[d] \ar[rd]&\\
		\pi_{\delbar}^{p,q}(X,x)\ar[rd]\ar[r,Rightarrow]&(\pi_{p+q}(X,x)\otimes\C)^\vee\ar[d]&\pi_{\del}^{p,q}(X,x)\ar[ld]\ar[l,Rightarrow]\\
		&\pi_A^{p,q}(X,x)&
	\end{tikzcd}
\end{equation}

\subsection{Holomorphic higher operations and formality}
As one might expect from the previous section, there are higher operations in complex geometry, the definition of which uses complex analytic information, e.g. that of solutions to the $\partial\bar\partial$-equation. The simplest instance are the Aeppli-Bott-Chern (ABC) triple Massey products, defined in \cite{AnToBCform}: For any three classes $[\alpha],[\beta],[\gamma]\in H_{BC}(X)$, s.t. $\alpha\beta=\del\delbar x$ and $\beta\gamma=\del\delbar y$, one defines the triple ABC-Massey product as
\[
\langle [\alpha],[\beta],[\gamma]\rangle_{ABC}\in H_A(X)/(\alpha H_A(X)+H_A(X)\gamma)
\]
Previously, Christopher Deninger has defined higher operations in real Deligne cohomology \cite{Den_higher}, which turn out the be closely related. Further operations with more than three inputs are produced in \cite{TaDiss}, \cite{MilSt_bigrform}. All these operations can be seen as invariants of the (real) bigraded quasi-isomorphism type of the cbba $A_X$ as is explained in \cite{MilSt_bigrform}. It should be noted that neither the ordinary rational homotopy type, nor the Dolbeault homotopy type of \cite{NT}, are sufficient to capture these kind of phenomena. In fact, in those theories any $\partial\bar\partial$-manifold is formal by \cite{DGMS}, \cite{NT}, while an important example of Sferruzza and Tomassini \cite{SfTo_DBC} shows that the triple ABC Massey products do not necessarily vanish on $\partial\bar\partial$ manifolds. This raised the following question

\begin{question}
	Are compact K\"ahler, or projective, manifolds pluripotentially formal? I.e. does there exist a chain of (real) bigraded-quasi isomorphisms of cbba's between $(A_X^{\ast,\ast},\del,\delbar,\wedge)$ and a bigraded algebra $H(X)$ with trivial differentials?\footnote{One could take $H=(H_{BC}^{\ast,\ast}(X),0,0,\wedge)$} 
\end{question}
One checks that pluripotential formality implies the vanishing of all higher operations such as ordinary, Dolbeault or ABC Massey products.

A positive answer to this question in some special cases (compact Hermitian symmetric spaces, K\"ahler manifolds with a Hodge diamond of complete intersection type) was given in \cite{MilSt_bigrform}, \cite{StePHT}. On the other hand, in very recent work \cite{PSZ24}, we show that in general, compact K\"ahler manifolds are quite far from being formal. This opens up a new homotopic toolbox for their study.

\begin{thm}[\cite{PSZ24}]\,
\begin{enumerate}
\item Any compact manifold with a surjective map to a Riemann surface of genus at least two supports a nontrivial ABC Massey product of the form $\langle \alpha,\alpha,\beta\rangle_{ABC}$, where $\alpha,\beta$ are $1$-forms. In particular this holds for any such curve $\Sigma_{g\geq 2}$.
	\item For any complex manifold of dimension $\geq 4$, there exists a finite sequence of blow-ups in points and lines such that the resulting manifold carries a nontrivial ABC Massey product of the form\linebreak $\langle D_1,D_2,D_3\rangle_{ABC}$, where the $D_i$ are divisor classes. The value of the Massey product (when paired with an appropriate class of complementary degree) is related to the cross ratio of four points on a line.
\end{enumerate}
\end{thm}

In the latter examples, the value of the Massey product varies if one holomorphically varies the configurations in which one blows up and this can be used to distinguish biholomorphism types of manifolds in families where all intermediate Jacobians, and the variation of Hodge structure given by the cohomology, are trivial.

Given how fruitful the study of the relative position of the Hodge decomposition on de Rham cohomology with respect to the rational structure given by singular cohomology is in K\"ahler and algebraic geometry, it seems not unreasonable to expect further applications from the interplay between the pluripotential homotopy type and the rational homotopy type.

\subsection*{Acknowledgments} This article is a modified and expanded version of the introduction to my habilitation thesis (with the same title). The ideas layed out here evolved over time and will likely evolve further. I benefited hugely from discussions with many people, including, but by no means limited to: D. Angella, J. Cirici, H. Kasuya, D. Kotschick, A. Milivojevic, G. Placini, D. Sullivan, S. Wilson, and L. Zoller. I thank the anonymous referees for useful comments which improved the presentation.


\noindent Jonas Stelzig, Mathematisches Institut der LMU M\"unchen, Theresienstra{\ss}e 39, 80333 M\"unchen, Germany. jonas.stelzig@math.lmu.de
\end{document}